\documentclass[11pt]{article}
\usepackage{amsthm, amsmath, amssymb, amsfonts, url, booktabs, tikz, setspace, fancyhdr, bm}

\usepackage{amsthm}

\usepackage{mathrsfs} 
\usepackage{xcolor}
\usepackage[margin = 2.3cm]{geometry}
\usepackage{hyperref, enumerate}
\usepackage[shortlabels]{enumitem}
\usepackage[babel]{microtype}
\usepackage[english]{babel}
\usepackage[capitalise]{cleveref}
\usepackage{comment}
\usepackage{bbm}
\usepackage{csquotes}
\usepackage{mathabx}
\usepackage{tikz}
\usetikzlibrary{arrows.meta,decorations.pathreplacing,positioning}
\usepackage{graphicx}
\usepackage{float}
\usepackage{thm-restate}
\usepackage{mathtools}

\counterwithin{figure}{section}

\newtheorem{conjecture}{Conjecture}

\newtheorem{theorem}{Theorem}[section]
\newtheorem*{thm-non}{Theorem}
\newtheorem{prop}[theorem]{Proposition}
\newtheorem{lemma}[theorem]{Lemma}

\newtheorem{claim}[theorem]{Claim}

\theoremstyle{definition}

\newtheorem*{defn-non}{Definition}

\newlist{Case}{enumerate}{2}
\setlist[Case, 1]{%
    label           =   {\bfseries Case \arabic*.},
    labelindent=1em ,labelwidth=1.3cm, labelsep*=1em, leftmargin =!
}
\setlist[Case, 2]{%
    label           =   {\bfseries Subcase \arabic{Casei}.\arabic*.},
    labelindent=-1em ,labelwidth=1.3cm, labelsep*=1em, leftmargin =!
}

\newenvironment{poc}{\begin{proof}[Proof of the claim]}{\end{proof}}

\usepackage{todonotes}

\title{Balanced Sperner families via the topological Tverberg theorem}
\author{
Chong Shangguan\thanks{Research Center for Mathematics and Interdisciplinary Sciences, Shandong University and Frontiers Science Center for Nonlinear Expectations, Ministry of Education, Qingdao, China. Email: theoreming@163.com.}
\and
Zixiang Xu\thanks{School of Mathematical Sciences, Zhejiang University, Hangzhou, China. Email: zixiangxu@zju.edu.cn.}
\and
Yulin Yang\thanks{Research Center for Mathematics and Interdisciplinary Sciences, Shandong University, Qingdao,
China. Email: forestyoung@mail.sdu.edu.cn.}
}
\date{}
\begin{document}
\maketitle

\begin{abstract}
For every prime power $r\ge 2$, we show that any Sperner family $\mathcal F\subseteq 2^{[n]}$ with $|\mathcal F|\ge (r-1)n+1$ contains $r$ pairwise disjoint nonempty subfamilies whose unions are all equal and whose intersections are all equal. 
For $r=2$, this confirms a conjecture of Heged\"{u}s, with the sharp threshold $n+1$. In this purely combinatorial problem, our proof combines a multilinear polynomial method, a continuity argument, and the topological Tverberg theorem.
\end{abstract}

\section{Introduction}
Let $[n]=\{1,\dots,n\}$. For $r\ge2$, a family $\mathcal{F}=\{F_1,\dots,F_m\}\subseteq 2^{[n]}$ is called \emph{$r$-balanced} 
if there are pairwise disjoint nonempty sets $I_1,\dots,I_r\subseteq[m]$ such that all unions $\bigcup_{i\in I_j}F_i$ are equal and all intersections $\bigcap_{i\in I_j}F_i$ are equal. 

Questions of this kind go back to work of Lindstr\"{o}m~\cite{lindstrom1972theorem} on balanced set systems. 
A classical observation says that every family $\mathcal F\subseteq 2^{[n]}$ of size at least $n+1$ is $\cup$-balanced, meaning that there exist two disjoint nonempty subfamilies that have the same union. 
Lindstr\"{o}m~\cite{lindstrom1993another} later proved a stronger result for $2$-balanced families: every family $\mathcal{F}\subseteq 2^{[n]}$ of size at least $n+2$ is $2$-balanced, where the threshold $n+2$ is best possible in general. 
Very recently,  Heged\"us~\cite{hegedus2026uniform} studied a uniform version of this theorem and proved that if $\mathcal F$ is uniform and $|\mathcal F|\ge n+1$, then $\mathcal F$ is $2$-balanced. 
Motivated by this improvement from $n+2$ to $n+1$, he proposed the following conjecture. 
Recall that $\mathcal{F}$ is \emph{Sperner} if no two distinct members are comparable under inclusion. 

\begin{conjecture}[\cite{hegedus2026uniform}]\label{conj:hegedus}
Let $\mathcal{F}\subseteq 2^{[n]}$ be a Sperner family. If $|\mathcal{F}|\ge n+1$,
then $\mathcal{F}$ is $2$-balanced.
\end{conjecture}
We prove Conjecture~\ref{conj:hegedus} in the following stronger form. 
The case $r=2$ is exactly Heged{\"u}s' conjecture.

\begin{theorem}\label{thm:main}
Let $r\ge2$ be a prime power and $n\ge1$ be an integer. Every Sperner family
$\mathcal{F}\subseteq2^{[n]}$ with $|\mathcal{F}|\ge (r-1)n+1$ is $r$-balanced.
\end{theorem}

Let us explain the idea behind the proof. 
Identify a set $F\subseteq[n]$ with its characteristic vector $\boldsymbol{1}_F\in\{0,1\}^n$. 
If two disjoint subfamilies have convex hulls that meet, then, after deleting the sets with zero coefficient in two convex representations of a common point, the remaining subfamilies are balanced. 
Indeed, in each coordinate, the common convex combination records whether that coordinate is absent from all participating sets, present in all of them, or present only in some of them. 
Thus Radon's theorem~\cite{radon1921mengen} immediately gives Lindstr\"om's bound $n+2$ for arbitrary families~\cite{lindstrom1993another}, 
where Radon's theorem states that any $n+2$ points in $\mathbb R^n$ have two disjoint subcollections with intersecting convex hulls. 
This also explains why the uniform case is easier. If all sets in $\mathcal F$ have size $k$, then all characteristic vectors lie in the affine hyperplane $\sum_{a=1}^n x_a=k$, which has dimension $n-1$. 
Radon's theorem in this hyperplane gives the threshold $n+1$, and hence recovers Heged\"us' uniform theorem~\cite{hegedus2026uniform} from this geometric viewpoint.

For general Sperner families, however, this direct affine dimension argument breaks down. 
One particular reason is that there exists some Sperner family of size $n+1$ such that the corresponding vectors span a space with full affine dimension $n$. 
Thus the main difficulty of Conjecture~\ref{conj:hegedus} is to find another $(n-1)$-dimensional structure behind an arbitrary Sperner family.

Our proof introduces such an $(n-1)$-dimensional structure, but it is not linear. 
Given an arbitrary Sperner family $\mathcal{F}\subseteq 2^{[n]}$ of size $n+1$, using the Sperner condition, 
we build a strictly coordinate-monotone function $g:[0,1]^n\rightarrow\mathbb{R}$ whose zero set contains all characteristic vectors $\boldsymbol{1}_{F}$ for $F\in\mathcal{F}$. 
This zero set then plays the role of a nonlinear hypersurface containing those vectors. 
We then construct a continuous map that sends points of the simplex into this hypersurface in a way that remembers their supports coordinatewise.
After embedding the hypersurface into an $(n-1)$-dimensional hyperplane, we apply a topological version of Radon's theorem~\cite{bajmoczy1979common} to settle Conjecture~\ref{conj:hegedus}. 
The same idea, with the topological Tverberg theorem~\cite{barany1981topological} in place of Radon's theorem, gives Theorem~\ref{thm:main} for every prime power $r$. 

Our result fits into a familiar theme: discrete problems often become clearer after one passes through continuous mathematics, 
as in Lov\'asz's proof of Kneser's conjecture via Borsuk--Ulam~\cite{borsuk1933drei,lovasz1978kneser} and Alon's necklace splitting theorem~\cite{alon1987splitting}. 
In this sense, the key ingredient of our proof is to apply a continuous argument to a purely combinatorial problem.

\section{Proof of Theorem~\ref{thm:main}}
We use the following standard form of the topological Tverberg theorem. 
For a positive integer $d$, let $\Delta^d$ denote the $d$-dimensional simplex, which is the convex hull of $d+1$ affinely independent points in $\mathbb{R}^d$. 
For a simplex $\Delta^d={\rm conv}\{\boldsymbol{v}_0,\boldsymbol{v}_1,\ldots,\boldsymbol{v}_d\}$ and a subset $I\subseteq\{0,1,\ldots,d\}$, 
the set $\Delta_I={\rm conv}\{\boldsymbol{v}_i:i\in I\}$ is called a {\it face} of $\Delta^d$, where ${\rm conv}$ denotes the convex hull of a set of points. 

\begin{theorem}[\cite{barany1981topological,Ozaydin1987equivariant,volovikov1996topological}]\label{thm:topological-tverberg}
Let $r$ be a prime power and $d\ge0$. For every continuous map $f:\Delta^{(r-1)(d+1)}\to\mathbb R^d$, 
there are $r$ pairwise disjoint nonempty faces $\Delta_{I_1},\dots,\Delta_{I_r}$ whose images have a common point.
\end{theorem}

It suffices to consider the case $|\mathcal{F}|=(r-1)n+1$, since any larger Sperner family contains a subfamily of this size. Put $N=(r-1)n$ and write $\mathcal{F}=\{F_0,F_1,\ldots,F_N\}$.

For a nonempty $S\subseteq\{0,1,\dots,N\}$, write $U(S)=\bigcup_{i\in S}F_i$ and $L(S)=\bigcap_{i\in S}F_i$, respectively. 
Define a box $C(S)\subseteq[0,1]^n$ by declaring its $a$-th coordinate interval to be
\begin{equation*}
C(S)_a=
\begin{cases}
\{0\}, & a\notin U(S),\\
\{1\}, & a\in L(S),\\
[0,1], & a\in U(S)\setminus L(S).
\end{cases}
\end{equation*}
This box records exactly which elements are in the union, which elements are in the intersection, and which elements occur only sometimes. 
We also need the following open part of this box:
\begin{equation*}
C^\circ(S)=\prod_{a=1}^n C^\circ(S)_a,
\end{equation*}
where
\begin{equation*}
C^\circ(S)_a=
\begin{cases}
\{0\}, & a\notin U(S),\\
\{1\}, & a\in L(S),\\
(0,1), & a\in U(S)\setminus L(S).
\end{cases}
\end{equation*}
Notice that $C^\circ(S)$ is obtained from $C(S)$ by replacing every closed interval $[0,1]$ by the open interval $(0,1)$, and leaving $\{0\}$ and $\{1\}$ unchanged.

\begin{lemma}\label{lem:box-equality}
Let $I_1,\dots,I_r$ be nonempty subsets of $\{0,1,\dots,N\}$. If
$\bigcap_{j=1}^r C^\circ(I_j)\ne\varnothing$, then all boxes $C(I_j)$
are equal. In particular, all unions $U(I_j)$ are equal and all
intersections $L(I_j)$ are equal.
\end{lemma}

\begin{proof}[Proof of Lemma~\ref{lem:box-equality}]
For each coordinate $a$, $C^\circ(S)_a$ is one of $\{0\}$, $\{1\}$, and $(0,1)$, and these three sets are disjoint. 
Thus a common point in the sets $C^\circ(I_j)$ forces the same type in every coordinate for every box. 
The last assertion is just the definition of $C(S)$.
\end{proof}

One of the key ingredients in this paper is the following construction. 
Informally, the Sperner condition allows us to place all characteristic vectors $\boldsymbol{1}_{F_i}$ on one strictly monotone level set. 
Formally, put $\varepsilon=\frac{1}{n+1}$, and let $w:2^{[n]}\to\mathbb R$ be
\begin{equation*}
w(A):=
\begin{cases}
0, & \textup{if }A\in\mathcal{F},\\[0.2em]
\varepsilon |A|, & \textup{if }A\text{ is incomparable with every member of }\mathcal{F},\\[0.2em]
\min\{|A|-|F|:\ F\in\mathcal{F},\ A\subsetneq F\}, & \text{if }A\subsetneq F\textup{ for some }F\in\mathcal{F},\\[0.5em]
\max\{|A|-|F|:\ F\in\mathcal{F},\ F\subsetneq A\}, & \textup{if }F\subsetneq A\textup{ for some }F\in\mathcal{F}.
\end{cases}
\end{equation*}
We then have the following simple consequence from the construction.

\begin{claim}\label{lem:weight-function}
If $A\subsetneq B$, then $w(A)<w(B)$.
\end{claim}

\begin{poc}
Fix $A\subsetneq B$.
Suppose first that $A\subsetneq F$ for some $F\in\mathcal F$. Then $w(A)<0$. 
If $B$ is not properly contained in any member of $\mathcal F$, then $w(B)\ge 0$, so $w(A)<w(B)$. 
If $B\subsetneq F$ for some $F\in\mathcal F$, then using $|A|<|B|$ it is easy to verify $w(A)<w(B)$.

If $A\in\mathcal F$, then $B$ properly contains a member of $\mathcal F$, namely $A$. Hence $w(B)\ge |B|-|A|>0=w(A)$.

If $A$ is incomparable with every member of $\mathcal F$, then $w(A)=\varepsilon |A|<1$. 
Moreover, $B$ cannot belong to, nor be properly contained in, a member of $\mathcal F$, for otherwise $A$ would be comparable with such a member. 
Thus either $B$ is also incomparable with every member of $\mathcal F$, giving $w(B)=\varepsilon |B|>\varepsilon |A|=w(A)$, or $B$ properly contains a member of $\mathcal F$, giving $w(B)\ge 1>w(A)$.

Finally, suppose that $F\subsetneq A$ for some $F\in\mathcal F$. 
Then $B$ also properly contains a member of $\mathcal F$. 
Again using $|A|<|B|$ it is easy to verify $w(A)<w(B)$.
Therefore $w(A)<w(B)$ in all cases.
\end{poc}

Now let $g:[0,1]^n\to\mathbb R$ be the multilinear extension of $w$:
\begin{equation*}
g(\boldsymbol{x})=\sum_{A\subseteq[n]}w(A)\prod_{a\in A}x_a\prod_{a\notin A}(1-x_a).
\end{equation*}
Then $g(\boldsymbol{1}_A)=w(A)$ for every $A\subseteq[n]$, 
so in particular $g(\boldsymbol{1}_{F_i})=0$ for all $i$. 
Moreover, for each $a\in [n]$,
\begin{equation*}
\frac{\partial g}{\partial x_a}(\boldsymbol{x})
=
\sum_{A\subseteq[n]\setminus\{a\}}
\bigl(w(A\cup\{a\})-w(A)\bigr)
\prod_{b\in A}x_b
\prod_{\substack{b\notin A\\ b\ne a}}(1-x_b).
\end{equation*}
Fix a coordinate $a$. For $A\subseteq[n]\setminus\{a\}$, write
\begin{equation*}
M_A(\boldsymbol{x})
:=
\prod_{b\in A}x_b
\prod_{\substack{b\notin A\\ b\ne a}}(1-x_b).
\end{equation*}
For $\boldsymbol{x}\in[0,1]^n$, each $M_A(\boldsymbol{x})$ is nonnegative. 
Also, as $A$ runs over all subsets of $[n]\setminus\{a\}$, 
these terms correspond exactly to choosing, for each $b\ne a$, either the factor $x_b$ or the factor $1-x_b$. 
Therefore
\begin{equation*}
\sum_{A\subseteq[n]\setminus\{a\}} M_A(\boldsymbol{x})
=
\prod_{b\ne a}(x_b+(1-x_b))
=1.
\end{equation*}
In particular, the numbers $M_A(x)$ are nonnegative and not all zero. Since $w(A\cup\{a\})-w(A)$ is positive by Claim~\ref{lem:weight-function}, 
we get $\frac{\partial g}{\partial x_a}(\boldsymbol{x})>0$ for every $x\in[0,1]^n$. 
It follows that increasing any one coordinate, while keeping the others fixed, strictly increases the value of $g$.

Set $Z=g^{-1}(0)$. We have $\boldsymbol{1}_{F_0},\boldsymbol{1}_{F_1},\dots,\boldsymbol{1}_{F_N}\in Z$. 
The next step is the main geometric construction of the proof. 
We want to replace the usual affine map\footnote{Here we mean that, in the uniform case $\mathcal{F}=\{F_0,F_1,\ldots,F_N\}\subseteq\binom{[n]}{k}$, it suffices to apply Theorem~\ref{thm:topological-tverberg} to the affine map $f:\Delta^N=\operatorname{conv}\{\boldsymbol{v}_0,\boldsymbol{v}_1,\ldots,\boldsymbol{v}_N\}\to\{\boldsymbol{x}\in\mathbb{R}^n:\sum_{a=1}^nx_a=k\}\cong\mathbb{R}^{n-1}$ defined by $f(\boldsymbol{v}_i)=\boldsymbol{1}_{F_i}$ for each $0\le i\le N$.}
from the simplex to the convex hull of the points $\boldsymbol{1}_{F_i}$ by a map whose image stays on the zero set $Z$, but which still remembers exactly which vertices are used. 
For a point $\boldsymbol{\lambda}\in\Delta^N$, written as $\boldsymbol{\lambda}=(\lambda_0,\lambda_1,\ldots,\lambda_N)$ in barycentric coordinates, denote its \emph{support} by $S_{\boldsymbol{\lambda}}:=\{i:\lambda_i>0\}$.
We shall construct a continuous map $q:\Delta^N\to Z$ such that
\begin{equation*}
q(\boldsymbol{\lambda})\in Z\cap C^\circ(S_{\boldsymbol{\lambda}})
\end{equation*}
for every $\boldsymbol{\lambda}\in\Delta^N$.

Let $\boldsymbol{z}_{\boldsymbol{\lambda}}=\sum_{i\in S_{\boldsymbol{\lambda}}}\lambda_i\boldsymbol{1}_{F_i}$.
Then, for each coordinate $a$, $(\boldsymbol{z}_{\boldsymbol{\lambda}})_a=\sum_{i\in S_{\boldsymbol{\lambda}}:\,a\in F_i}\lambda_i$.
If $S_{\boldsymbol{\lambda}}=\{i\}$, set $q(\boldsymbol{\lambda})=\boldsymbol{1}_{F_i}$. 
Now assume $|S_{\boldsymbol{\lambda}}|\ge2$. 
The parameter $\theta>0$ will move the free coordinates (that is, the coordinates labeled by $U(S_{\boldsymbol{\lambda}})\setminus L(S_{\boldsymbol{\lambda}})$) monotonically from $0$ to $1$. 
Define $\boldsymbol{\gamma}_{\boldsymbol{\lambda}}(\theta)\in[0,1]^n$ by
\begin{equation*}
\boldsymbol{\gamma}_{\boldsymbol{\lambda}}(\theta)_a
=
\begin{cases}
0, & a\notin U(S_{\boldsymbol{\lambda}}),\\
1, & a\in L(S_{\boldsymbol{\lambda}}),\\
\dfrac{\theta\cdot(\boldsymbol{z}_{\boldsymbol{\lambda}})_a}
{1-(\boldsymbol{z}_{\boldsymbol{\lambda}})_a+\theta\cdot(\boldsymbol{z}_{\boldsymbol{\lambda}})_a}, &
a\in U(S_{\boldsymbol{\lambda}})\setminus L(S_{\boldsymbol{\lambda}}).
\end{cases}
\end{equation*}

In the third case $a\in U(S_{\boldsymbol{\lambda}})\setminus L(S_{\boldsymbol{\lambda}})$, 
it is easy to see that $a$ belongs to at least one of the sets $F_i$, $i\in S_{\boldsymbol{\lambda}}$, but not to all of them. 
Since all $\lambda_i$ with $i\in S_{\boldsymbol{\lambda}}$ are positive and sum to $1$, this gives
\begin{equation*}
0<(\boldsymbol{z}_{\boldsymbol{\lambda}})_a=\sum_{i\in S_{\boldsymbol{\lambda}}:\,a\in F_i}\lambda_i<1.
\end{equation*}
Hence $\boldsymbol{\gamma}_{\boldsymbol{\lambda}}(\theta)_a\in(0,1)$ for every $\theta>0$, and so
$\boldsymbol{\gamma}_{\boldsymbol{\lambda}}(\theta)\in C^\circ(S_{\boldsymbol{\lambda}})$.

As $\theta\to0$, all free coordinates tend to $0$, while the coordinates forced by $L(S_{\boldsymbol{\lambda}})$ remain equal to $1$. 
Therefore $\boldsymbol{\gamma}_{\boldsymbol{\lambda}}(\theta)\to\boldsymbol{1}_{L(S_{\boldsymbol{\lambda}})}$. 
Similarly, as $\theta\to\infty$, all free coordinates tend to $1$, and hence $\boldsymbol{\gamma}_{\boldsymbol{\lambda}}(\theta)\to\boldsymbol{1}_{U(S_{\boldsymbol{\lambda}})}$.

Moreover, since $\mathcal F$ is Sperner and $|S_{\boldsymbol{\lambda}}|\ge2$, every $i\in S_{\boldsymbol{\lambda}}$ satisfies
\begin{equation*}
L(S_{\boldsymbol{\lambda}})\subsetneq F_i\subsetneq U(S_{\boldsymbol{\lambda}}).
\end{equation*}
Indeed, if $L(S_{\boldsymbol{\lambda}})=F_i$, then $F_i\subseteq F_j$ for every $j\in S_{\boldsymbol{\lambda}}$, contradicting the Sperner property. 
Similarly, if $F_i=U(S_{\boldsymbol{\lambda}})$, then $F_j\subseteq F_i$ for every $j\in S_{\boldsymbol{\lambda}}$, again a contradiction. 
Therefore, by Claim~\ref{lem:weight-function} and the strict coordinate monotonicity of $g$, we obtain
\begin{equation*}
g(\boldsymbol{1}_{L(S_{\boldsymbol{\lambda}})})
<g(\boldsymbol{1}_{F_i})=0
<g(\boldsymbol{1}_{U(S_{\boldsymbol{\lambda}})}).
\end{equation*}
On the other hand, each free coordinate of $\boldsymbol{\gamma}_{\boldsymbol{\lambda}}(\theta)$ is a strictly increasing function of $\theta$, while the fixed coordinates stay fixed. 
Therefore $\theta\mapsto g(\boldsymbol{\gamma}_{\boldsymbol{\lambda}}(\theta))$ is strictly increasing. 
Since its limits at $0$ and $\infty$ have opposite signs, it has exactly one zero. 
We denote this zero by $\theta({\boldsymbol{\lambda}})$, and set
\begin{equation*}
q({\boldsymbol{\lambda}})=\boldsymbol{\gamma}_{\boldsymbol{\lambda}}(\theta({\boldsymbol{\lambda}})).
\end{equation*}
We then show that the map $q$ constructed above is continuous. 
The only issue is that the support $S_{\boldsymbol{\lambda}}$ may shrink when some coordinates of ${\boldsymbol{\lambda}}$ tend to $0$.

\begin{claim}\label{claim:q-continuous}
The map $q$ is continuous.
\end{claim}

\begin{poc}
Let $\boldsymbol{\lambda}^{(k)}\to\boldsymbol{\lambda}$ in $\Delta^N$.
We show that $q(\boldsymbol{\lambda}^{(k)})\to q(\boldsymbol{\lambda})$.
Since $q(\boldsymbol{\lambda}^{(k)})\in [0,1]^n$ and $[0,1]^n$ is compact, 
it suffices to show that every subsequential limit of $q(\boldsymbol{\lambda}^{(k)})$ is equal to $q(\boldsymbol{\lambda})$. 
Thus, after passing to an arbitrary convergent subsequence (which we do not relabel), we may assume that $q(\boldsymbol{\lambda}^{(k)})\to \boldsymbol y$. 
We shall prove that $\boldsymbol y=q(\boldsymbol{\lambda})$. 
Since $Z=g^{-1}(\{0\})$ is the preimage of the closed set $\{0\}$ under the continuous function $g$, the set $Z$ is closed. Hence $\boldsymbol y\in Z$.

Moreover, the supports $S_{\boldsymbol{\lambda}^{(k)}}$ can take only finitely many values. 
Passing to a further subsequence (which we do not relabel), we may assume that $S_{\boldsymbol{\lambda}^{(k)}}=T$ is constant. 
Put $S=S_{\boldsymbol{\lambda}}$. Then $S\subseteq T$. 
If $|T|=1$, then all $\boldsymbol{\lambda}^{(k)}$ are the same vertex of $\Delta^N$, and the claim is clear. 
We may therefore assume that $|T|\ge2$.

Write $\theta_k=\theta(\boldsymbol{\lambda}^{(k)})$. Then each $\theta_k$ is a positive real number. 
By compactness of the extended interval $[0,\infty]$, after passing, if necessary, to a further subsequence of the one just chosen (which we do not relabel), 
we may assume that $\theta_k\to \tau$ for some $\tau\in[0,\infty]$.

For convenience, we record the coordinatewise formula of $\boldsymbol{y}$.
\begin{equation*}
\boldsymbol{y}_a=\lim_{k\to\infty}\boldsymbol{\gamma}_{\boldsymbol{\lambda}^{(k)}}(\theta_k)_a
=\lim_{k\to\infty}
\begin{cases}
0, & a\notin U(T), \\
1, & a\in L(T), \\
\dfrac{\theta_k\cdot(\boldsymbol{z}_{\boldsymbol{\lambda}^{(k)}})_a}{1-(\boldsymbol{z}_{\boldsymbol{\lambda}^{(k)}})_a+\theta_k\cdot (\boldsymbol{z}_{\boldsymbol{\lambda}^{(k)}})_a}, & a\in U(T)\setminus L(T).
\end{cases}
\end{equation*}

First, suppose $0<\tau<\infty$. 
Since $\boldsymbol{z}_{\boldsymbol{\lambda}^{(k)}}\to\boldsymbol{z}_{\boldsymbol{\lambda}}$, 
we have $(\boldsymbol{z}_{\boldsymbol{\lambda}^{(k)}})_a\to(\boldsymbol{z}_{\boldsymbol{\lambda}})_a=0$ for $a\notin U(S)$, 
and $(\boldsymbol{z}_{\boldsymbol{\lambda}^{(k)}})_a\to (\boldsymbol{z}_{\boldsymbol{\lambda}})_a=1$ for $a\in L(S)$.
Therefore, taking the limit coordinatewise gives
\begin{equation*}
\boldsymbol{y}_a=
\begin{cases}
0, & a\notin U(S), \\
1, & a\in L(S), \\
\dfrac{\tau\cdot(\boldsymbol{z}_{\boldsymbol{\lambda}})_a}{1-(\boldsymbol{z}_{\boldsymbol{\lambda}})_a+\tau\cdot (\boldsymbol{z}_{\boldsymbol{\lambda}})_a}, & a\in U(S)\setminus L(S).
\end{cases}
\end{equation*}
If $|S|\ge2$, this gives $\boldsymbol{y}=\boldsymbol{\gamma}_{\boldsymbol{\lambda}}(\tau)$.
As $\boldsymbol{y}\in Z$, we have $g(\boldsymbol{\gamma}_{\boldsymbol{\lambda}}(\tau))=0$. 
By the uniqueness of the zero of $\theta\mapsto g(\boldsymbol{\gamma}_{\boldsymbol{\lambda}}(\theta))$, 
it follows that $\tau=\theta(\boldsymbol{\lambda})$, and hence $\boldsymbol{y}=q(\boldsymbol{\lambda})$. 
If $S=\{i\}$, this gives $\boldsymbol{y}=\boldsymbol{1}_{F_i}=q(\boldsymbol{\lambda})$.

Next, suppose $\tau=0$. We claim that $\boldsymbol{y}\le \boldsymbol{1}_{L(S)}$ coordinatewise. 
Indeed, let $a\notin L(S)$. Then $(\boldsymbol{z}_{\boldsymbol{\lambda}})_a<1$. 
If $a\notin U(T)$, then the corresponding coordinate is constantly $0$. 
Otherwise, since $L(T)\subseteq L(S)$, we have $a\in U(T)\setminus L(T)$, and 
\begin{equation*}
\frac{\theta_k\cdot(\boldsymbol{z}_{\boldsymbol{\lambda}^{(k)}})_a}
{1-(\boldsymbol{z}_{\boldsymbol{\lambda}^{(k)}})_a
+\theta_k\cdot(\boldsymbol{z}_{\boldsymbol{\lambda}^{(k)}})_a}
\to0.
\end{equation*}
This proves the claim. 
If $|S|\ge 2$, then $g(\boldsymbol{y})\le g(\boldsymbol{1}_{L(S)})<0$, contradicting $\boldsymbol{y}\in Z$. 
Hence $S=\{i\}$ for some $i$. Then $\boldsymbol{y}\le \boldsymbol{1}_{F_i}$ and $g(\boldsymbol{y})=g(\boldsymbol{1}_{F_i})=0$. 
Thus, by the strict monotonicity of $g$, we must have $\boldsymbol{y}=\boldsymbol{1}_{F_i}=q(\boldsymbol{\lambda})$.

Finally, suppose $\tau=\infty$. We claim that $\boldsymbol{y}\ge \boldsymbol{1}_{U(S)}$ coordinatewise. 
Indeed, let $a\in U(S)$. Then $(\boldsymbol{z}_{\boldsymbol{\lambda}})_a>0$. 
If $a\in L(T)$, then the corresponding coordinate is constantly $1$. 
Otherwise, since $U(S)\subseteq U(T)$, we have $a\in U(T)\setminus L(T)$, and 
\begin{equation*}
\frac{\theta_k\cdot(\boldsymbol{z}_{\boldsymbol{\lambda}^{(k)}})_a}
{1-(\boldsymbol{z}_{\boldsymbol{\lambda}^{(k)}})_a
+\theta_k\cdot(\boldsymbol{z}_{\boldsymbol{\lambda}^{(k)}})_a}
\to1.
\end{equation*}
This proves the claim. 
If $|S|\ge2$, then $g(\boldsymbol{y})\ge g(\boldsymbol{1}_{U(S)})>0$, again contradicting $\boldsymbol{y}\in Z$. 
Hence $S=\{i\}$ for some $i$. Then $\boldsymbol{y}\ge\boldsymbol{1}_{F_i}$ and $g(\boldsymbol{y})=g(\boldsymbol{1}_{F_i})=0$. 
Therefore, the strict monotonicity gives $\boldsymbol{y}=\boldsymbol{1}_{F_i}=q(\boldsymbol{\lambda})$.

Thus every subsequential limit of $q(\boldsymbol{\lambda}^{(k)})$ is $q(\boldsymbol{\lambda})$. 
By compactness of $[0,1]^n$, it follows that $q(\boldsymbol{\lambda}^{(k)})\to q(\boldsymbol{\lambda})$. Hence $q$ is continuous.
\end{poc}

Let $H=\left\{\boldsymbol{x}\in\mathbb R^n:\sum_{a=1}^n x_a=0\right\}$.
We project $Z$ to $H$ along the direction $\boldsymbol{1}=(1,\dots,1)$.
More precisely, define
\begin{equation*}
\pi(\boldsymbol{x})
=
\boldsymbol{x}
-\frac1n\left(\sum_{a=1}^n x_a\right)\boldsymbol{1}.
\end{equation*}

\begin{claim}\label{claim:PiInj}
    $\pi$ is injective on $Z$.
\end{claim}

\begin{poc}
Suppose that $\boldsymbol{x},\boldsymbol{y}\in Z$ and $\pi(\boldsymbol{x})=\pi(\boldsymbol{y})$. 
Then there exists $c\in\mathbb{R}$ such that $\boldsymbol{x}=\boldsymbol{y}+c\boldsymbol{1}$.
Interchanging $\boldsymbol{x}$ and $\boldsymbol{y}$, if necessary, we may assume that $c\ge 0$. 
Hence $\boldsymbol{x}\ge \boldsymbol{y}$ coordinatewise. 
Since $g(\boldsymbol{x})=g(\boldsymbol{y})=0$, 
the strict monotonicity of $g$ implies that $\boldsymbol{x}$ cannot be strictly larger than $\boldsymbol{y}$. 
Therefore $c=0$, and hence $\boldsymbol{x}=\boldsymbol{y}$.
\end{poc}

We next compose the map $q:\Delta^N\to Z$ with the projection $\pi$ onto $H$ and then obtain a continuous map
\begin{equation*}
\pi\circ q:\Delta^N\rightarrow H.
\end{equation*}
Since $H$ has dimension $n-1$ and $N=(r-1)n$, 
we apply Theorem~\ref{thm:topological-tverberg} with $d=n-1$. 
Thus there are pairwise disjoint nonempty faces
\begin{equation*}
\Delta_{J_1},\dots,\Delta_{J_r}\subseteq\Delta^N
\end{equation*}
and points $\boldsymbol{\lambda}^{(a)}\in\Delta_{J_a}$, $a=1,\dots,r$, such that
\begin{equation*}
\pi(q(\boldsymbol{\lambda}^{(1)}))
=
\pi(q(\boldsymbol{\lambda}^{(2)}))
=
\cdots
=
\pi(q(\boldsymbol{\lambda}^{(r)})).
\end{equation*}
Since all the points $q(\boldsymbol{\lambda}^{(a)})$ lie in $Z$ and $\pi$ is injective on $Z$, we have
\begin{equation*}
q(\boldsymbol{\lambda}^{(1)})
=
q(\boldsymbol{\lambda}^{(2)})
=
\cdots
=
q(\boldsymbol{\lambda}^{(r)}).
\end{equation*}
For each $a=1,\dots,r$, put $I_a=S_{\boldsymbol{\lambda}^{(a)}}$.
Since the faces $\Delta_{J_1},\dots,\Delta_{J_r}$ are pairwise disjoint, the sets $I_1,\dots,I_r$ are pairwise disjoint and nonempty. 
By the defining property of $q$, $q(\boldsymbol{\lambda}^{(a)})\in C^\circ(I_a)$ for every $a=1,2,\ldots,r$.
But the points $q(\boldsymbol{\lambda}^{(a)})$ are all equal, so $\bigcap_{a=1}^r C^\circ(I_a)\ne\varnothing$.
Lemma~\ref{lem:box-equality} now gives
\begin{equation*}
U(I_1)=U(I_2)=\cdots=U(I_r)
\end{equation*}
and
\begin{equation*}
L(I_1)=L(I_2)=\cdots=L(I_r).
\end{equation*}
This finishes the proof.

\section{Some remarks}
In this paper we prove Conjecture~\ref{conj:hegedus}. 
In particular, the threshold $n+1$ in Heged\"us' original conjecture is sharp. 
For example, one can easily check that $\mathcal{F}=\bigl\{\{1\},\{2\},\dots,\{n\}\bigr\}$ is Sperner and not $2$-balanced.

We also remark that the prime power assumption in Theorem~\ref{thm:main} comes from our use of the topological Tverberg theorem, 
and this statement is known to fail in general when $r$ is not a prime power, see~\cite{blagojevic2019barycenters,frick2015counterexamples,skopenkov2018user}. 
Thus our method does not directly extend to arbitrary $r$. 
It remains natural, however, to ask whether the combinatorial bound in Theorem~\ref{thm:main} is close to best possible. 
The following construction does not quite reach $(r-1)n$, but it comes within a constant depending only on $r$, and gives some evidence for this expectation.
\begin{prop}\label{prop:improved-construction}
For every $r\ge 3$ and every $n\ge 2\lceil\log_2(r-1)\rceil$, 
there exists a Sperner family $\mathcal G\subseteq 2^{[n]}$ which is not $r$-balanced and satisfies
\begin{equation*}
|\mathcal G|\ge (r-1)n-2(r-1)\lceil\log_2(r-1)\rceil .
\end{equation*}
\end{prop}
\begin{proof}[Proof of Proposition~\ref{prop:improved-construction}]
Put $s=\lceil\log_2(r-1)\rceil$. 
For each $i\in [n]$, we consider the $(s+1)$-sets whose largest element is exactly $i$, that is, the sets $A\cup\{i\}$ with $A\in \binom{[i-1]}{s}$. 
From all such $(s+1)$-sets, we arbitrarily pick $\min\left\{r-1,\binom{i-1}{s}\right\}$ sets to form a family $\mathcal{G}$. 
Clearly $\mathcal{G}$ is $(s+1)$-uniform, hence $\mathcal{G}$ is Sperner.

Since $n\ge 2s$, we have $\binom{2s}{s}\ge 2^s\ge r-1$. 
Hence for every $i\ge 2s+1$, there are at least $r-1$ available $(s+1)$-sets with largest element $i$. 
Therefore, from each such $i$, our construction chooses exactly $r-1$ sets. 
It follows that
\begin{equation*}
|\mathcal G|
=\sum_{i=1}^n \min\left\{r-1,\binom{i-1}{s}\right\}
\ge (n-2s)(r-1)
\ge (r-1)n-2(r-1)\lceil\log_2(r-1)\rceil .
\end{equation*}

It remains to show that $\mathcal{G}$ is not $r$-balanced. 
Suppose that there are pairwise disjoint nonempty subfamilies $\mathcal H_1,\dots,\mathcal H_r\subseteq\mathcal G$ 
whose unions are all equal and whose intersections are all equal. 
Let $i$ be the largest element appearing in some member of $\mathcal H_1\cup\cdots\cup\mathcal H_r$, and denote the common union by $U$. 
Since $i\in U$, each $\mathcal H_j$ must contain at least one set containing $i$. 
As the subfamilies $\mathcal H_1,\dots,\mathcal H_r$ are pairwise disjoint, 
this gives at least $r$ distinct members of $\mathcal G$ containing $i$. 
On the other hand, by the choice of $i$, every member of $\mathcal H_1\cup\cdots\cup\mathcal H_r$ which contains $i$ has largest element exactly $i$. 
But in the construction we chose at most $r-1$ members of $\mathcal G$ whose largest element is $i$, which is a contradiction. 
This finishes the proof.
\end{proof}

\section*{Acknowledgements}
Zixiang Xu is grateful to Chong Shangguan for inviting him to visit Shandong University in May 2026 and for his hospitality during the visit.
Chong Shangguan and Yulin Yang are supported by the National Natural Science Foundation of China under Grant Nos. 12571352 and 12231014, 
and the Fundamental Research Funds for the Central Universities.

\bibliographystyle{abbrv}
\bibliography{ref}
\end{document}